\documentclass[12pt,a4paper]{article}
\usepackage{fullpage,makeidx,amsfonts,amssymb,amstex,epic}
\newcommand{\C}{{\mathbb C}}
\newcommand{\Q}{{\mathbb Q}}
\newcommand{\N}{{\mathbb N}}

\renewcommand{\P}{{\mathbb P}}

\newcommand{\R}{{\mathbb R}}

\newcommand{\ka}{{\cal A}}
\newcommand{\kb}{{\cal B}}

\newcommand{\kf}{{\cal F}}

\newcommand{\kh}{{\cal H}}
\newcommand{\ki}{{\cal I}}
\newcommand{\kj}{{\cal J}}

\newcommand{\kl}{{\cal L}}
\newcommand{\km}{{\cal M}}
\newcommand{\kn}{{\cal N}}
\newcommand{\ko}{{\cal O}}
\newcommand{\kp}{{\cal P}}
\newcommand{\kq}{{\cal Q}}

\newcommand{\kv}{{\cal V}}

\begin{document}
\newtheorem{lemma}{Lemma}[section]
\newtheorem{proposition}[lemma]{Proposition}
\newtheorem{remark}[lemma]{Remark}
\newtheorem{example}[lemma]{Example}
\newtheorem{theorem}[lemma]{Theorem}
\newtheorem{definition}[lemma]{Definition}
\newtheorem{corollary}[lemma]{Corollary}


\title{{\Large\bf  Semicontinuity for representations of one--dimensional
    Cohen--Macaulay Rings}}
\author{
Y.A. Drozd*\\
Mechanics \&\\
Mathematics Faculty\\
Kiew University\\
Wladimirskaya 64\\
252017 Kiew\\
Ukraine
\and
G.-M. Greuel**\\
Fachbereich Mathematik\\
Universit\"at Kaiserslautern\\
Erwin-Schr\"odinger-Stra\ss e\\
67663 Kaiserslautern\\
Federal Republic of Germany}

\maketitle
\thispagestyle{empty}
\vspace{0.5cm}
\tableofcontents
\vspace{1cm}

$\ast$ Partially supported by the Foundation for Fundamental
Research of the Ukraine and the Deutsche Forschungsgemeinschaft (DFG)\\
$\ast\ast$ Partially supported by the Deutsche Forschungsgemeinschaft (DFG)\\
\thispagestyle{empty}
\setcounter{page}{0}



\section*{Introduction}\addcontentsline{toc}{section}{Introduction}

``Algebraic families'' of modules and algebras play an important role
in several questions of representation theory.  It is often especially
useful to know that some ``discrete invariants'' are constant or, at
least, are semi--continuous in such families, that is they can change
only in ``exceptional points'' which form a family of smaller
dimension.  Perhaps the best known results in this direction are those
of Gabriel \cite{Gab} and Kn\"orrer \cite{Kn}.  Gabriel proved that
finite representation type is an open condition for finite dimensional
algebras (``fat points''), while Kn\"orrer showed that the number of
parameters for modules of prescribed rank is semi--continuous in
families of commutative Cohen--Macaulay rings of Krull dimension 1
(``curve singularities'').  In \cite{DG 2} Kn\"orrer's theorem was
used to show that the unimodal singularities of type $T_{pq}$ are of
tame Cohen--Macaulay type.

Unfortunately, the arguments of \cite{Kn} do not work in the
non--commutative case.  The aim of this paper is to refine them in
such a way that they could be applied to non--commutative
Cohen--Macaulay algebras, too.  For this purpose we introduce the
notion of ``dense subrings'' which seems rather technical but,
nevertheless, useful.  It enables the construction of ``almost
versal'' families of modules for a given algebra (cf.\ Theorem
\ref{3.4}) and the definition of the ``number of parameters''.  Just
as in the commutative case, it is important that the bases of these
``almost versal'' families are projective varieties.  Once having
this, we are able to prove an analogue of Kn\"orrer's theorem (cf.\
Theorem \ref{4.9}) and a certain variant (cf.\ Theorem 4.11) which
turns out to be useful, for instance, to extend the tameness criterion
for commutative algebras \cite{DG 2} to the case of characteristic 2.
The semicontinuity implies, in particular, that the set of so--called
``wild algebras'' in any family is a countable union of closed
subsets.  A very exciting problem is whether it is actually closed,
hence whether the set of tame algebras is open.  However, Theorem 4.9,
together with the results of \cite{DG 2}, imply that tame is indeed an
open property for curve singularities (commutative one-dimensional
Cohen--Macaulay rings).

An analogous procedure leads to the semicontinuity of the number of parameters
in other cases, like representations of finite dimensional algebras or
elements  of finite dimensional bimodules.

Though we do not consider here the problem of constructing moduli
spaces for Cohen-Macaulay modules (cf. \cite{GP}), we may rephrase the
semicontinuity theorem by saying that the dimension of the moduli
space for such modules of prescribed rank varies upper semicontinuosly
in flat families of Cohen-Macaulay algebras. Likewise, the
semicontinuity in other cases (finite-dimensional algebras or
bimodules) may be also understood  as semicontinuity of the dimension
of the corresponding moduli spaces of representations under
deformations of the algebra or bimodule.

Unfortunately, our results are, just as all known results till now, not
sufficient to prove the ``tame is open condition''--conjecture. Nevertheless,
the semicontinuity theorem as well as the construction of ``almost versal''
families of modules are not restricted to tame algebras. They have a
potentially broader field of applications to classification problems in
representation theory. They are a particularly powerful tool if, for a given
algebra, the deformation theory of this algebra is sufficiently known and the
classification problem for the deformed algebras is easier to solve or even
known. The great success of this approach in the commutative case is, of
course, also due to the fact that the deformation theory of singularities is a
highly developed field. We hope that this paper stimulates further research in
the deformation theory for non--commutative Cohen--Macaulay algebras.

 \newpage

\section{Dense subalgebras}

\begin{definition}\label{1.1}
Let $B$ be a subring of a ring $A$.   Call $B$ {\bf dense}
in $A$ if any simple $A$--module $U$ is also simple as $B$--$End_A
U$--bimodule.
\end{definition}

{\bf Examples}:
\begin{enumerate}
\item If $A$ is commutative, then $End_A U = A/Ann\, U$, so any subring $B$ is
dense in $A$.

\item We can take $A = M_2 (\R)$ the ring of $2 \times 2$ real
matrices and $B = \C$ (or even $B = \Q(i))$, naturally embedded in is
$A$.  More generally, let $k$ be some field, $L, K$ and $B$ be its
extensions such that $L = KB$ (a composite) and $(L : K) = n$.  Then
we have a natural embedding $L \to A = M_n(K)$, thus also $B$ is a
subalgebra of $A$ and $B$ is dense in $A$.

\end{enumerate}

\begin{lemma}\label{1.2}
Let $D$ be a division algebra over an infinite field $k$, $A = M_n(D)$ and $B
\subset A$ a dense subalgebra.
Consider $W = M_{n \times m}(D)$ as $A$--module and let $V \subset W$ be
a $B$--submodule such that $AV = W$.   Suppose that $m = nq + r$ with $0 \le r
< n$.   Then there exists an automorphism $\sigma$ of $W$ such that $\sigma(V)$
contains the matrices
\begin{eqnarray*}
E_1 & = & (I\; 0 \ldots 0\; 0'),\;\; E_2 = (0\; I \ldots 0\; 0'), \ldots,\\
E_q & = & (0\; 0 \ldots 0\; I\; 0')
\end{eqnarray*}
and a matrix of the form $(Y_1\; Y_2 \ldots Y_q\; Y')$.

\begin{tabular}{lll}
Here & $I$ & denotes the $n \times n$ unit matrix,\\
     & $0$ & denotes the $n \times n$ zero matrix, \\
     & $0'$& denotes the $n \times r$ zero matrix,
\end{tabular}

$Y_1, Y_2, \ldots, Y_q$ are some $n \times n$ matrices and $Y'$ is an $n \times
r$ matrix of rank $r$ (of course, if $r = 0$, then $0'$ and $Y'$ are empty, so
in this case $V$ contains an $A$--basis of $W \simeq q A)$.
\end{lemma}

{\bf Proof:}  Use induction on $m$.   First prove the claim in the case $m
\le n$.   Choose a matrix $X \in V$ of maximal possible rank, say $d$.
Then we must show that $d = m$.   Suppose that $d < m$.   Denote by $\bar{X}_1,
\bar{X}_2, \ldots, \bar{X}_m$ the columns of $X$ and let $\bar{X}_1,
\bar{X}_2, \ldots, \bar{X}_d$ be linear independent.   Then there exists an
automorphism $\sigma$ of $W$ such that the last $(m-d)$ columns of $\sigma(X)$
are zero.   So we may suppose that $\bar{X}_{d+1} = \ldots = \bar{X}_m =
\bar{0}$.   Note that $W \simeq mU$ where $U = nD$ is the only simple
$A$--module.   As $d< n,\; \bar{X}_1, \bar{X}_2, \ldots, \bar{X}_d$ do not
span $U$ over $D$.   Denote by $V'$ the projection of $V$ onto the $(d +
1)$--st component of $W$ (i.e.\ $V' \subset U$ consists of the first $d + 1$
columns of the matrix from $V$).   As $AV = W$, we have $AV' = U$, so $V'
\not= 0$.   But as $U$ is a simple $B$--$D$--bimodule, $V'D = U$.   Therefore,
$V$ contains a matrix $Y$ such that its $(d + 1)$--st column $\bar{Y}_{d+1}$
does not lie in $\langle\bar{X}_1,\; \bar{X}_2, \ldots, \bar{X}_d \rangle D$.
Choose a
$D$--basis of $U$ of the form $\{\bar{X}_1,\; \bar{X}_2, \ldots, \bar{X}_d, \;
\bar{Y}_{d+1}, \bar{Z}_{d+2}, \ldots, \bar{Z}_n\}$ and let
\[
\bar{Y}_i = \sum^d_{j=1} \bar{X}_i \lambda_{ij} + \bar{Y}_{d+1} \lambda_{i,
d+1} + \sum^n_{j=d+2} \bar{Z}_j \lambda_{ij}.
\]

Again, using an automorphism of $W$, we may suppose that $\lambda_{i, d+1} = 0$
for $i \not= d + 1$.  But then rank$(\gamma X + Y) \ge d + 1$ for some $\gamma
\in k$, which is a contradiction.
Hence, $d = m$.   In particular, if $m = n$, there exists an automorphism
$\sigma$ of $W$, such that $\sigma(X) = I$.   Thus, our claim is proved for $m
\le n$.

Suppose now that $m > n$ and consider the projection $V'$ of $V$ onto $W' =
nU$, the first $n$ components of $W$ (that is the first $n$  columns of each
matrix $X \in W)$.   As we have proved, there exists an automorphism $\sigma'$
of $W'$ such that $\sigma'(V') \ni I$.   We can extend $\sigma'$ to $W$ and
thus suppose that $V$ contains a matrix $X$ of the form ($I
X'$).   Again using an automorphism, we obtain that $X' = 0$, that is $X =
E_1$.   Now consider the projection $V''$ of $V$ onto $W'' = (m-n)U$, the last
components of $W$.   Using induction, we may also suppose that the claim is
valid for $V''$, thus $V$ contains the matrices of the form:
\[
\begin{array}{ll}
(X_2 \;\; I \;\; 0 \ldots 0 \;\; 0'),\quad & (X_3\;\; 0\;\; I \ldots 0\;\;
0'),\; \ldots,\\
(X_q\:\; 0 \ldots 0 \;\; I\;\; 0'),\quad & (Y_1\;\; Y_2 \ldots Y_q\;\;
Y')
\end{array}
\]
with $rank(Y') = r$.   But then, again using an automorphism of $W$, we can
make $X_2 = X_3 = \ldots = X_q = 0$, q.e.d.

\begin{corollary}\label{1.3}
Let $B \subset A$ be a dense subring.   Suppose that $A/rad\, A$ is an
artinian ring containing an infinite field $k$ in its centre and, moreover, $k
\subset B/(B \cap rad\, A)$ (for example $B$ and $A$ are $k$--algebras).   Let
$V \subset nA$ be a $B$--submodule such that $AV = nA$.   Then $V$ contains an
$A$--basis of $nA$.
\end{corollary}

{\bf Proof:}  Of course, we may replace $A$ by $A/rad\, A$, so suppose that $A
= \prod^s_{i=1} A_i$ with $A_i$ simple artinian.  Put $B_i = pr_i\, B,\; V_i =
pr_i\, V$, $pr_i$ being the projection from $A$ onto $A_i$.   Then $B_i$ is
dense in $A_i$ and $A_iV_i = nA_i$.   Lemma
\ref{1.2} implies that each $V_i$ contains an $A_i$--basis $\{\bar{e}_{ij}|
j= 1, \ldots, n\}$ of $nA_i$.   Let $e_{ij} \in V$ be such elements that
$pr_ie_{ij} = \bar{e}_{ij}\; (j = 1, \ldots, n;\;\; i = 1, \ldots, s)$.
Consider in $V$ the elements $e_j(\lambda_1, \ldots, \lambda_s) =
\sum^s_{i=1} \lambda_i e_{ij}$ where $\lambda_1,  \ldots, \lambda_s
\in k$.   The sets $T_i = \{(\lambda_1,  \ldots, \lambda_s) | pr_i
e_j(\lambda_1, \ldots, \lambda_s)$ form a basis of $nA_i\}$, are Zariski--open
in $k^s$ and non--empty.   As $k$ is infinite, their intersection is also
non--empty.

But if $(\lambda_1, \ldots, \lambda_s)$ lies in this intersection,
then $\{e_j (\lambda_1, \ldots, \lambda_s)|j = 1, 2, \ldots, n\}$
is a basis of $nA$, q.e.d.

\begin{definition}\label{1.4}

Let $D$ be a skewfield, $A = M_n(D),\;\; U = nD$ the simple $A$--module and
\[
\kf:\; U = U_0 \supset U_1 \supset \ldots \supset U_s = \{0\}
\]
a flag of $D$--subspaces in $U$.   Put $A(\kf) = \{a \in A \mid a U_i \subset
U_i$ for all $i = 0, 1, \ldots, m\}$ and call $A(\kf)$ a {\bf flag subalgebra}
in $A$.   If $A = \prod_i A_i$ with $A_i$ simple artinian, call a {\bf flag
subalgebra} of $A$ any subalgebra $A'$ of the form $A' = \Pi_i A'_i$ where
$A'_i$ is a flag subalgebra in $A_i$ for each $i$.
\end{definition}

\begin{lemma}\label{1.5}

Let $A$ be a semi--simple artinian ring,
$A(\kf)$ a
flag subalgebra in $A$ and $B$ a subring of
$A(\kf)$.
Then there exists a flag subalgebra $A'$ such that $B \subset A' \subset
A(\kf)$ and $B$ is dense in $A'$.
\end{lemma}

{\bf Proof}:  Obviously, we may suppose $A$ to be simple.   Then take a maximal
$B$--invariant flag $\kf' \supset
\kf$ and put $A' = A(\kf')$.

\begin{proposition}\label{1.6}

Let $A$ be an algebra over a separably  closed field $k$, such that $A/rad\,A$
is finite--dimensional over $k$, $B$ a dense subalgebra of $A$ and $K$ a
separably
generated extension of $k$.   Then $B \otimes_k K$ is also dense in $A
\otimes_k K$.
\end{proposition}

{\bf Proof:}  Of course, we may suppose $A$ to be simple finite--dimensional,
that is $A = M_n(F)$ for some skewfield $F$.   As $k$ has no separable
extensions, $F$ is really a field \cite{DK}, so it is a pure inseparable
extension of $k$.   But then $F \otimes_k K$ is again a field, hence $A
\otimes_k K$ is simple.   Its only simple module is $U \otimes_k K$, where $U$
is the simple $A$--module and $F \otimes_k K$ is its endormorphism ring.   But
the same observation shows that for any simple $B$--$F$--bimodule $V$ (e.g.\
for $U$) the tensor product $V \otimes_k K$ remains simple, q.e.d.
\newpage



\section{Cohen--Macaulay Algebras}

We consider here one--dimensional Cohen--Macaulay algebras (not necessarily
commutative), also known as orders in semi--simple algebras.
\begin{definition}\label{2.1}

Call a ring $\Lambda$ a {\bf CM--algebra} (more precisely, {\rm{\bf1}}--{\bf
  dimensional, analytically reduced Cohen--Macaulay algebra}) if it satisfies
  the following conditions:

\begin{enumerate}
\item [(1)] $\Lambda$ is an algebra over a one--dimensional local, commutative,
noetherian ring $R$, which is
a finitely generated and torsion--free $R$--module.

\item [(2)] The completion $\widehat{\Lambda}$ of $\Lambda$ in the ${\frak
    m}$--adic topology, where ${\frak  m}$ is the maximal ideal of $R$,
    contains no nilpotent ideals.
\end{enumerate}
\end{definition}

It follows from (2) that, in this case, $R$ is a Cohen--Macaulay ring and
$\Lambda$ is a maximal Cohen--Macaulay $R$--module.

In particular, the ${\frak  m}$--adic completion $\widehat{R}$ of $R$ has no
nilpotent elements.   Denote by $Q$ (respectively $\widehat{Q}$) the total ring
of fractions of $R$ (respectively $\widehat{R}$).   Then both $Q$ and
$\widehat{Q}$ are finite products of fields and $Q\Lambda = Q \otimes_R
\Lambda\; (\widehat{Q} \Lambda = \widehat{Q} \otimes_R \Lambda = \widehat{Q}
\otimes_{\widehat{R}} \widehat{\Lambda}$) is a semi--simple artinian
$Q$--algebra (respectively $\widehat{Q}$--algebra).   If $\Gamma$ is a subring
of $Q\Lambda$, containing $\Lambda$ and being finitely generated as
$\Lambda$--module (or, equivalently, as $R$--module), call it an {\bf overring}
of $\Lambda$.   Of course, any such overring is also a CM--algebra.
If $\Lambda$ has no proper overrings, call it a {\bf maximal} CM--algebra.
It is known (cf., e.g., \cite{D1}) that, under condition (1),  condition
(2) is equivalent to the existence of maximal overrings of $\Lambda$.   More
precisely, under conditions (1) and (2), the overrings of $\Lambda$ satisfy
the ascending chain conditions and any two maximal overrings of $\Lambda$ are
conjugate in $Q\Lambda$.

Let $\Lambda$ be a CM--algebra.  We call a $\Lambda$--module $M$ a
$\Lambda$--{\bf lattice} (or a {\bf Cohen--Macaulay}{\boldmath$\Lambda$}{\bf
  --module}) provided it is a maximal Cohen--Macaulay $R$--module.  Denote
by CM($\Lambda$) the category of all $\Lambda$--lattices.  Any such lattice $M$
embeds naturally into the finitely generated $Q\Lambda$--module $QM =
Q\otimes_RM$.  So, if $\Gamma$ is an overring of $\Lambda$, the
$\Gamma$--module $\Gamma M \subset QM$ is well--defined.

The following assertions are  rather well--known (for the case when $R$ is
a discrete valuation ring, cf. \cite{Rog};  the proofs in the general situation
are  the same).

\begin{proposition}\label{2.2}

\begin{itemize}
\item[(a)] Any maximal CM--algebra $\Lambda$ is hereditary (that is
gl.$\dim\,\Lambda = 1$ or, equivalently, any $\Lambda$--lattice is projective).

\item[(b)] Let $\Lambda$ be a maximal CM--algebra and $A$ any flag subalgebra
of $\Lambda/rad\,\Lambda$ (cf.\ Definition \ref{1.4}).   Then the preimage
of $A$ in $\Lambda$ is hereditary and any hereditary CM--algebra can be
obtained in this way.
\end{itemize}
\end{proposition}

\begin{corollary}\label{2.3}

Let $\Omega$ be a hereditary (e.g.\ maximal) overring of a CM--algebra
$\Lambda$.   Then there exists a hereditary overring $\Omega'$ such that
$\Lambda \subset \Omega' \subset \Omega$ and $\Lambda$ is dense in $\Omega'$.
\end{corollary}

\begin{proposition}\label{2.4}

Suppose that the residue field $k = R/{\frak  m}$ is infinite.   Let $\Gamma$
be an
overring of $\Lambda$ such that $\Lambda$ is dense in $\Gamma$ and $M$ be a
Cohen--Macaulay $\Lambda$--module such that $\Gamma M \simeq n \Gamma$.   Then
$M$ is isomorphic to a module $M'$ such that $n \Lambda \subset M' \subset n
\Gamma$.
\end{proposition}

{\bf Proof:}  We may suppose that $M \subset n \Gamma$.   By Corollary
\ref{1.3} it contains a basis of $n \Gamma$.   Then there exists an
automorphism $\sigma$ of $n\Gamma$ which maps this basis to the standard one,
namely $(1, 0, \ldots, 0),\; (0, 1, 0, \ldots, 0), \ldots, (0,
\ldots, 0, 1)$.   Therefore, $M' = \sigma(M) \supset n \Lambda$.

\begin{proposition}\label{2.5}

If $M, M' \subset n\Gamma$ are $\Lambda$--submodules such that $\Gamma M =
\Gamma M' = n\Gamma$ (e.g.\ $M$ and $M'$ contain $n\Lambda$), then $M \simeq
M'$ if and only if there exists an automorphism $\sigma \in Aut(n\Gamma)$ such
that $\sigma(M) = M'$.
\end{proposition}

The proof is evident.

\begin{definition}\label{2.6}

For any Cohen--Macaulay $\Lambda$--module $M$ denote $\ell(M)$ the length of
the $Q\Lambda$--module $QM$ and call it the {\bf rational length} of $M$.
\end{definition}

{\bf Remark:} If $\Gamma$ is an overring of $\Lambda$ and $M$ is a
Cohen--Macaulay $\Gamma$--module, the rational length of $M$ does not depend
on whether we consider $M$ as a $\Lambda$-- or as a $\Gamma$--module.   On the
other hand, we have to distinguish between $\ell(M)$ and $\ell(\widehat{M})$
where $\widehat{M}$ is the ${\frak  m}$--adic completion of $M$.

Recall some connections between Cohen--Macaulay modules and their completions.
The proofs can be found in \cite{CR} or \cite{Rog} for the case when
$R$ is a discrete valuation ring, and they are also valid in the general
situation.

\begin{proposition}\label{2.7}

\begin{enumerate}
\item[(a)] $M \simeq N$ if and only if $\widehat{M} \simeq \widehat{N}$.
\item[(b)] If $N$ is a Cohen--Macaulay $\widehat{\Lambda}$--module such that
  $QN \simeq Q\widehat{N}'$ for some Cohen--Macaulay $\Lambda$--module $N'$,
  then there exists a Cohen--Macaulay $\Lambda$--module $N''$ such that $N
  \simeq \widehat{N}''$.
\item[(c)] If $\widehat{N}$ is isomorphic to a direct summand of $\widehat{M}$,
  then $N$ is isomorphic to a direct summand of $M$.
\end{enumerate}
\end{proposition}

In the next section we shall use the following simple result.

\begin{proposition}\label{2.8}
Let $P$ be a projective $\Lambda$--module.   Then there exists a projective
$\Lambda$--module $P'$ such that $P \oplus P'$ is free of rank $r \le
\dim_k(P/rad\, P)$ where $k = R/{\frak  m}$.
\end{proposition}

{\bf Proof:}  Due to Proposition \ref{2.7}, we may suppose that $R$ is
complete, thus the Krull--Schmidt theorem holds for modules.   Let $\Lambda
\simeq \oplus^s_{i=1} n_i P_i$, where all $P_i$ are indecomposable and
pairwise non--isomorphic.   Then $P\simeq \oplus^s_{i=1} m_i P_i$ for some
$m_i$.

Take $r$ the least integer such that $r n_i \ge m_i$ for all $i$.   Then
$r\Lambda \simeq P \oplus P'$ for $P' = \oplus^s_{i=1}(rn_i - m_i)P_i$.   As
$\dim_k(P/rad\,P) = \sum^s_{i=1} m_i \dim_k(P_i/rad\, P_i) \ge m_i$, one
obtains $r \le \dim_k (P/rad\, P)$, q.e.d.

{\bf Remark:}  Obviously, $\dim_k(P/rad\,P) \le \ell(\widehat{P})$, so the last
number can also serve as an upper bound for $r$.
\newpage


\section{Families of Modules}

{}From now on we suppose that our rings are algebras over the field $k =
R/{\frak  m}$, where $R$ is, as in the preceding paragraph, a
Cohen-Macaulay ring.

\begin{definition}\label{3.1}
Let $X$ be a $k$--scheme, $\ko_X = \ko$ its structure sheaf, $\Lambda$ a
CM--algebra (1--dimensional and analytically reduced) and $\km$ a coherent
sheaf on $X$ of $\Lambda \otimes_k \ko$--modules.   Call $\km$ a {\bf family}
of Cohen--Macaulay $\Lambda$--modules on $X$ if the following conditions hold:

\begin{enumerate}
\item[(1)] $\km$ is $R$--torsion free.
\item[(2)]$\km$ is $\ko$--flat.
\item[(3)] For each point $x \in X,\; \km(x) = \km \otimes_{\ko}
k(x)$ is a Cohen--Macaulay $\Lambda(x)$--module, where $\Lambda(x) = \Lambda
\otimes_k k(x)$.
\end{enumerate}

It is easy to see that, under conditions (1) and (2), condition (3) is
equivalent to:

\begin{enumerate}
\item[(3')] For every non--zero divisor $a \in R$, the sheaf
$\km/a\km$ is also $\ko$--flat.
\end{enumerate}
\end{definition}

We are going to construct some ``almost universal'' families.   Let $\Gamma$
be an overring of $\Lambda$ (cf.\ \S 2) and fix some positive integers $n$ and
$d$.   Put $\Phi = \Gamma/\Lambda$ and consider the Grassmannian $Gr =
Gr(n\Phi,d)$, that is
the variety of subspaces of codimension $d$ in $n\Phi$.   Recall that for every
$k$--scheme $X$ the morphisms $X \to Gr$ are in 1--1 correspondence with
$\ko$--factormodules of $n\Phi \otimes_k \ko_X$ which are locally free of
rank $d$ \cite{Mum}.   Consider the subvariety $B = B(n, d; \Lambda, \Gamma)$
of $Gr(n\Phi, d)$ consisting of all $\Lambda$--submodules of $n\Phi$.   In
other words, the morphisms $X \to B$ are in 1--1 correspondence with $\Lambda
\otimes_k \ko_X$--factormodules of $n\Phi\otimes_k \ko_X$ which are locally
free over $\ko_X$ of rank $d$.   Evidently it is a closed subscheme of $Gr$.
Denote by $\kf = \kf(n, d; \Lambda, \Gamma)$ the preimage in $n\Gamma
\otimes_k \ko_B$ of the canonical locally free sheaf of corank $d$ on $B$.  As
$(n\Gamma\otimes_k\ko_B)/\kf$ is flat over $\ko_B$, one can see that $\kf$ is
really a family of Cohen--Macaulay $\Lambda$--modules on $B$ having the
following universal property (cf.\ \cite{GP}).

\begin{proposition}\label{3.2}
For any family of Cohen--Macaulay $\Lambda$--modules $\km$ on a scheme $X$
such that $n\Lambda \otimes_k \ko_X \subset \km \subset n\Gamma \otimes_k
\ko_X$ and $(n\Gamma \otimes_k \ko_X)/\km$ is locally free over $\ko_X$ of rank
$d$, there exists a unique morphism $\varphi : X \to B$ such that $\km =
\varphi^\ast(\kf)$.
\end{proposition}

\begin{definition}\label{3.3}
Call the families satisfying the conditions of Proposition \ref{3.2}  {\bf
sandwiched families} with respect to $\Gamma$ of rank $n$ and codimension
$d$.   In particular, when $X =$ Spec $k$, we have {\bf sandwiched
modules} (with respect to $\Gamma$).
\end{definition}

{}From now on we suppose the ground field $k$ to be algebraically closed.   We
are going to show that the sandwiched families are, in some sense ``almost
versal'', that is any other families can be stably  glued from finitely many
sandwiched families.   Taking into account Corollary \ref{2.3}, this follows
from the following result.

\begin{theorem}\label{3.4a}
Let $\Gamma$ be a hereditary overring of $\Lambda$ such that $\Lambda$ is
dense in $\Gamma$.   Then, given a family $\km$ of Cohen--Macaulay
$\Lambda$--modules on a reduced $k$--scheme $X$, there exists a descending
chain of closed subschemes $X = X_0 \supset X_1 \supset X_2 \supset \cdots
\supset X_m = \emptyset$, a set of morphisms $\{\varphi_i : Y_i
\longrightarrow B(n_i, d_i;\, \Lambda, \Gamma) \,|\,i = 1, \ldots, m\}$ and a
set of projective  $\Gamma$--modules $\{P_i\,|\,i = 1, \ldots, m\}$ such that
$\km_{Y_i} \oplus (P_i \otimes_k \ko_{Y_i}) \simeq \varphi^\ast_i \kf(n_i,
d_i;\, \Delta, \Gamma)$ where $Y_i = X_{i-1} \backslash X_i$ and $n_i \le
\hat{\ell}(\km) = \ell(\widehat{\km}(x))$ for an arbitrary closed point $x \in
X$.
\end{theorem}

Indeed, we shall establish a more general result, when $\Gamma$ is not
necessarily hereditary, but $\Gamma\km$ is flat over $\Gamma\otimes_k\ko_X$.
If $\Gamma$ is hereditary, the last condition becomes superfluous.   Since
$\Gamma\km/\km$ is $\ko_X$--coherent, there exists an open dense subset $U
\subseteq X$, on which $\Gamma\km/\km$ is flat over $\ko$.   Then, as $\Gamma$
is hereditary, it follows from \cite{CE} (Theorem \ref{2.8}) that $\Gamma\km$
is also $\Gamma\otimes_k\ko$--flat.   Moreover, the function $x \mapsto
\dim_{k(x)} \Gamma\km(x)/\mbox{rad}\Gamma\km(x)$ takes its maximum in some
closed point of $X$ and it does not exceed $\hat{\ell}(\km)$.   Hence, we need
only to establish the following fact:

\begin{proposition}\label{3.4}
Let $\Gamma$ be an overring of $\Lambda$ such that $\Lambda$ is dense in
$\Gamma$.      Let $\km$ be a family of
Cohen--Macaulay $\Lambda$--modules on a reduced $k$--scheme $X$ such that
$\Gamma\km$ is flat over
$\Gamma \otimes_k\ko$.   Then there exists an open subscheme $Y \subset X$, a
projective $\Gamma$--module $P$ and a morphism $\varphi : Y \to B (n, d;
\Lambda, \Gamma)$ for some integers $n$ and $d$ such that the restriction on
$Y$ of the family $\km \oplus (P \otimes_k \ko)$ is isomorphic to
$\varphi^\ast \kf(n, d; \Lambda, \Gamma$).   Moreover, we can choose $n \le
\max_g \ell(\km(g))$ where $g$ runs through minimal points of $X$ (that is
generic points of its irreducible components).
\end{proposition}

{\bf Proof:}  Of course, we may suppose that $X$ is irreducible.   Let $g \in
X$ be its generic point.   Consider the $\Gamma$--module $\Gamma \km(g)$.   It
is finitely generated and flat over $\Gamma(g)$, hence projective \cite{Bou},
(Ch.\ I, \S 2, Ex.\ 15).
By Proposition \ref{2.8}, there exists a projective
$\Gamma(g)$--module $P'$ such that $\Gamma \km(g) \oplus P' \simeq n
\Gamma(g)$ and we can choose $n
\le \dim_{k(g)} (\Gamma \km (g)/rad\, \Gamma \km(g))$.   If we move  to the
completions, there is a 1--1--correspondence between projective and
semi--simple $\widehat{\Gamma}$--modules and the same is valid for
$\widehat{\Gamma}(g)$--modules.   But as $k$ is separably closed and $k(g)$
separably generated over $k$, we have seen in the proof of Proposition
\ref{1.6} that any simple $\widehat{\Gamma}(g)$--module is of the form $U
\otimes_k k(g)$ for some simple $\widehat{\Gamma}$--module $U$.   Hence, the
same is true for projectives, so $\widehat{P}' \simeq \widehat{P} \otimes_k
k(g)$ for some projective $\widehat{\Gamma}$--module $\widehat{P}$.   But
Proposition \ref{2.7}
implies then that $\widehat{P}$ is really a completion of some projective
$\Gamma$--module $P$, whence $P' \simeq P \otimes_k k(g)$.

Replacing $\km$ by $\km \oplus (P \otimes_k \ko)$, we may now suppose that
$\Gamma \km(g) \simeq n \Gamma(g)$.   But $\Lambda(g)$ is dense in $\Gamma(g)$
by Proposition \ref{1.6}, so we may suppose, using  Corollary \ref{1.3}, that
$\km(g)$ contains a  basis of $n \Gamma(g)$.
By Proposition \ref{2.4}, $\km(g)$ is isomorphic to a submodule of $n
\Gamma(g)$ containing $n \Lambda(g)$.   So let $n\Lambda (g) \subset \km(g)
\subset n \Gamma(g)$.   Then the same is true on an open subset $Y \subset X$,
that is $n \Lambda \otimes \ko_Y \subset \km_Y \subset n \Gamma \otimes
\ko_Y$.   Shrinking $Y$, we may also suppose that $(n \Gamma \otimes
\ko_Y)/\km_Y$ is locally free of some rank $d$ (over $\ko_Y$) and it remains
to use Proposition \ref{3.2}.

An obvious iteration gives us the necessary generalization of Theorem 3.4:

\begin{corollary}\label{3.5}
Under the conditions of Proposition \ref{3.4} there exists a descending chain
of closed subschemes $X = X_0 \supset X_1 \supset X_2 \supset \ldots \supset
X_n = \emptyset$, a set of morphisms $\{\varphi_i : Y_i \to B(n_i, d_i,
\Lambda, \Gamma) | i = 1, 2, \ldots, n\}$ and a set of projective
$\Gamma$--modules
$\{P_i | i = 1, 2, \ldots, n\}$ such that $\km_{Y_i} \oplus (P_i \otimes_k
\ko_{Y_i}) \simeq \varphi^\ast_i \kf (n_i, d_i; \Lambda, \Gamma)$ where $Y_i =
X_{i-1} \backslash X_i$ and $n_i \le \max \{\dim_{k(x)} (\Gamma \km(x)/rad\,
\Gamma \km(x)) | x \in X\}$.
\end{corollary}

Fix now a  CM--algebra $\Lambda$
and an  overring $\Gamma$.   Let $B = B(n, d;
\Lambda, \Gamma)$ and $\kf = \kf(n, d; \Lambda, \Gamma)$.  Choose a two--sided
$\Gamma$--ideal $I \subset rad\,\Lambda$ of finite codimension (over $k$) and
put $F = \Gamma/I;\; \bar{\Lambda} = \Lambda/I$.   We can identify $Gr(n\Phi,
d)$ with the closed subscheme of $Gr(nF, d)$ consisting of all
subspaces $V$ containing $n\bar{\Lambda}$.   Then $B$ also becomes a closed
subvariety of $Gr(nF, d)$.   We shall consider the elements of $nF$ as
rows of length $n$ with entries from $F$ and identify $Aut(nF)$ with the full
linear group $GL(n, F) = G$ acting on $nF$ according to the rule  $g \cdot v
= vg^{-1}$.   Then Proposition \ref{2.5} implies that two subspaces $V, V' \in
B$ correspond to isomorphic sandwiched modules if and only if there exists an
element $g \in G$ such that $g \cdot V = V'$.

Considering the elements of $nF$ as the rows of $n \times n$ matrices, we can
identify $nV$ with a subspace in $M_n(F)$.   Then we obtain the following:

\begin{proposition}\label{3.6}
Let $V \in B,\; g \in G$.   Then $g \cdot V \in B$ if and only if $g \in G
\cap nV$.   Hence, $G \cdot V \cap B = (G \cap nV) \cdot V \simeq (G \cap
nV)/StV$, where $StV = \{g \in G | gV = V\}$.
\end{proposition}

As $G$ is open in $M_n(F),\; G \cap nV$ is open in $nV$, hence $\dim(G\cap nV)
= \dim nV = n(\gamma n - d)$ where $\gamma = \dim\, F$.   Therefore, $\dim(G
\cdot V \cap B) = n(\gamma n - d)  - \dim\, StV$, whence:

\begin{corollary}\label{3.7}
For each integer $i$ the set $B_i = \{V \in B | \dim(G\cdot V \cap B) \le i\}$
is closed in $B$.
\end{corollary}

Put
\[
{\rm par}(n, d; \Lambda, \Gamma) = \max_i(\dim B_i - i)
\]
\mbox{ and}
\[ {\rm par}(n; \Lambda,
\Gamma)= \max_d {\rm par}(n, d; \Lambda, \Gamma).
\]
Intuitively, ${\rm par}(n, d;
\Lambda, \Gamma)$ is the number of independent parameters defining the
isomorphism classes of
sandwiched $\Lambda$--modules of rank $n$ and codimension $d$ with respect to
$\Gamma$.   Corollary \ref{3.5} evidently implies the following result.

\begin{corollary}\label{3.7a}
Under the conditions of Proposition 3.5, for any closed point $x \in X$ the
set $\{y \in X \mid \km(y) \simeq \km(x) \otimes_k k(y)\}$ is constructible
(that is a finite union of locally closed subsets of $X$) and its dimension is
bigger or equal to dim$X -$ {\rm par}$(\widehat{\ell}(\km);\, \Delta, \Gamma$).
\end{corollary}

In particular, this assertion is true for any family of Cohen--Macaulay
$\Lambda$--modules if we take for $\Gamma$ an hereditary overring of $\Lambda$
such that $\Lambda$ is dense in $\Gamma$ (which always exists, cf.\ Corollary
2.3).

\begin{corollary}\label{3.9}
Let $\Gamma$ be any overring of $\Lambda$ and $\Omega$ a hereditary
overring of $\Lambda$ such that $\Lambda$ is dense in $\Omega$.   Put $\ell_0
= \ell(\widehat{\Lambda})$.   Then ${\rm par}(n, \Gamma) \le {\rm par}(\ell_0n,
\Omega)$ for all $n$.
\end{corollary}

Of course, if $\Gamma \subset \Gamma'$ are two overrings of $\Lambda$ and
$\dim_k(\Gamma'/\Gamma)= c$, then ${\rm par}(n, d; \Lambda, \Gamma) \le {\rm
  par}(n, c + d; \Lambda, \Gamma')$, whence ${\rm par}(n; \Lambda, \Gamma) \le
{\rm par}(n; \Lambda, \Gamma')$.   Put
\[
b(n, \Lambda) = \max \{{\rm par}(n; \Lambda, \Gamma)\}
\]
where $\Gamma$ runs through all overrings of $\Lambda$ (we have actually to
look only for maximal ones).   Let also $p(n, \Lambda)$ denote the maximal
value of $\dim\, X - \dim \{y \in X \mid \km(y) \simeq \km(x) \otimes_k k(y)\}$
taken for all families $\km$ with all possible bases $X$ and for all closed
points $x \in X$.

\begin{corollary}\label{3.10}
Let $\ell_0 = \ell(\widehat{\Lambda})$.   Then
\[
b(n, \Lambda) \le p(n, \Lambda) \le b (n\ell_0, \Lambda).
\]
\end{corollary}
\newpage

\section{Families of algebras}

Now we formulate and prove the semicontinuity statements in two
variants: for ``familes of algebras'' (Theorems~\ref{4.7} and \ref{4.9}) and
for ``families of generators'' (Theorem~\ref{411}).

Again $k$ denotes an algebraically closed field.

\begin{definition}\label{4.1}
Let $C$ be a reduced algebraic curve over  $k$, $\Lambda$ a coherent sheaf of
$\ko_C$--algebras, containing no nilpotent ideals and such that for every
point $p \in C$, $\Lambda_p$ is maximal Cohen--Macaulay $\ko_{C,p}$--module.
Then we call $\Lambda$ a {\bf sheaf of CM--algebras} or just a {\bf CM--algebra
  on  $\bf C$}.
\end{definition}

If $\Gamma \supset \Lambda$ is another CM--algebra on $C$ and, for each $p \in
C$, $\Gamma_p$ is overring of $\Lambda_p$, call $\Gamma$ {\bf an  overring} of
$\Lambda$.

\begin{proposition}\label{4.2}
If $\Lambda$ is a CM--algebra on a curve $C$, then, for every point $p \in C$,
$\Lambda_p$ is a CM--algebra (in the sense of Definition \ref{2.1}) and,
moreover, for almost all points $\Lambda_p$ is maximal.
\end{proposition}

{\bf Proof:}  We only need to prove that $\widehat{\Lambda}_p$ contains no
nilpotent ideal.   According to \cite{D1}, this is equivalent to the existence
of a maximal overring of $\Lambda_p$.   Denote by $Z$ the centre of
$\Lambda_p$. It is a localization of a finitely generated $k$--algbra, hence,
its algebraic closure
$\bar{Z}$ in the total quotient ring $Q$ is a finitely generated $Z$--module
(cf.\ \cite{Bou}, Ch.V.\ \S 3.2).   As before, we consider $\Lambda_p$ embedded
in $Q \Lambda_p = Q
\otimes_Z \Lambda_p$.   Therefore, $\bar{Z}\Lambda_p \subset Q \Lambda_p$ is
well--defined.   But now $Q\Lambda_p$ is a central, semi--simple, hence,
separable $QZ$--algebra, so $\bar{Z}\Lambda_p$ has a maximal overring (cf.\
\cite{CR}).   Moreover, as $\bar{Z}_p = Z_p$ for almost all $p \in C$ and
$\bar{Z}\Lambda_p$ is maximal for almost all $p$, the same is true also for
$\Lambda_p$, q.e.d.

Call $\Lambda$ {\bf hereditary} if all $\Lambda_p$ are hereditary (note that,
for a general point $g \in C$, $\Lambda_g \simeq Q\Lambda$ is semi--simple).
It is well--known (cf.\ \cite{CR}) that one--dimensional CM--algebras can be
defined locally:

\begin{proposition}\label{4.3}
Let $\Lambda$ be a CM--algebra on a curve $C$ and suppose that for each closed
point $p \in C$ an overring $\Gamma (p) \supset \Lambda_p$ is given such that
$\Gamma(p) = \Lambda_p$ for almost all $p$.   Then there exists an
overring $\Gamma \supset \Lambda$ such that $\Gamma_p = \Gamma(p)$ for all $p$.
\end{proposition}

\begin{corollary}\label{4.4}
There exists a hereditary overring $\Omega \supset \Lambda$ such that
$\Lambda_p$ is dense in $\Omega_p$ for each $p \in C$.
\end{corollary}

Let now $\Gamma$ be any overring of $\Lambda$.   As $\Gamma_p = \Lambda_p$ for
almost all $p$, the sum
\[
{\rm par}(n; \Lambda, \Gamma) = \sum_{p \in C} {\rm par}(n; \Lambda_p,
\Gamma_p)
\]
is well--defined.

\begin{definition}\label{4.5}
Let $f : Y \to X$ be a morphism of $k$--schemes and $\kl$ be a coherent sheaf
of $\ko_Y$--algebras.   Call $(\kl, f) = (\kl, f : Y \to X)$ a {\bf family of
  CM--algebras} with the base $X$ provided the following conditions hold:

\begin{enumerate}
\item[(1)] $f$ is flat and $f_\ast(\kl)$ is flat $\ko_X$--module.
\item[(2)] $Y(x) = f^{-1}(x)$ is a reduced algebraic curve for each $x \in X$.
\item[(3)] $\kl(x)$ is a CM--algebra on $Y(x)$ for each $x \in X$.
\end{enumerate}
\end{definition}

\begin{definition}\label{4.6}
Let $(\kl, f : Y \to X)$ be a family of CM--algebras with base $X$.   A
{\bf family of overrings} of $(\kl, f)$ is a family $(\kl', f)$ (with the same
$f$) such that $\kl' \supset \kl,\; f_\ast(\kl'/\kl)$ is $\ko_X$--flat and, for
each $x \in X$,  $\kl'(x)$ is an overring of $\kl(x)$.
\end{definition}

Given a family of overrings $\kl' \supset \kl$, we can define the functions on
$X$:
\begin{eqnarray*}
{\rm par}(x,n,d) & := & {\rm par}(x,n,d; \kl, \kl') = {\rm par}(n, d; \kl(x),
\kl'(x)),\\
{\rm par}(x, n)&  := & {\rm par}(x, n; \kl, \kl') = {\rm par}(n; \kl(x),
\kl'(x)).
\end{eqnarray*}

\begin{theorem}\label{4.7}
The functions ${\rm par}(x, n, d)$ and ${\rm par}(x,n)$ are
upper-semicontinuous, that is for each integer $i$ and for any $k$--scheme $X$
the sets $X_i(d) = \{x \in X | {\rm par}(x, n, d) \ge i\}$ and
$X_i = \{x \in X|{\rm par}(x,n) \ge i\}$ are closed in $X$.
\end{theorem}

{\bf Proof:}  As $X_i = \cup_d X_i(d)$ and since this union is finite, we only
need to prove that $X_i(d)$ is closed.   Moreover, we may suppose that $X$ is
a smooth curve.   Let $\kn = \kl'/\kl$.   Consider the relative Grassmannian
$Gr(n\kn, d) \to X$ and its closed subscheme (over $X$) $\kb(n, d)$ consisting
of $\kl$--submodules.   Let $\kj$ be the biggest two--sided $\kl'$--ideal
contained in $\kl$.   Then it is easy to see that $\kl/\kj$ is torsion--free
over $\ko_X$, hence, flat.   Thus, $\kl'/\kj$ is also flat over $\ko_X$.   As
in the proof of Proposition \ref{3.6}, identify $Gr(n\kn, d)$ with the closed
subscheme of $Gr(n\bar{\kl}', d)$, where $\bar{\kl}' = \kl'/\kj$, and consider
the group scheme over $X$, $GL(n, \bar{\kl}')$ acting on the last Grassmannian.
The same observations as in the proof of Proposition \ref{3.6} shows that
$\kb_j = \{v \in \kb(n, d) \mid \dim\, St\, v \ge j\}$ is closed in $\kb$.   As
$\kb$ is proper over $X$, its projection $Z_j$ is also closed.   But, by
definition $X_i = \cup_j X_{ij}$, where $X_{ij} = \{x \in Z_j \mid \dim\,
\kb_j(x) \ge i + j\}$ are closed, q.e.d.

{\bf Remark:}  Suppose that the base of the family $(\kl, f)$ is a smooth
curve and both $\kl$ and $\kl'$ are Cohen--Macaulay $\ko_Y$--modules.   Then
$(\kl', f)$ is a family of overrings as it follows from \cite{BG} (Example
3.2.5).   Moreover, in this case $\ko_Y$ is Cohen--Macaulay itself and $\dim\,
Y = 2$.   Hence, we are able to construct Cohen--Macaulay $\ko_Y$--modules
locally as in the following lemma.

\begin{lemma}\label{4.8}
Suppose that $Y$ is a reduced 2--dimensional Cohen--Macaulay variety.   Let
$\km$ be a Cohen--Macaulay $\ko_Y$--module, $\{y_1, y_2, \ldots, y_m\}$ a set
of points of $Y$ of codimension 1 and $N(y_i)$ a finitely generated
$\ko_{Y,y_i}$--submodule in $\kq\km$ where $\kq$ is the total quotient ring of
$\ko$.   Then there exists the Cohen--Macaulay submodule $\kn \subset \kq\km$
such that $\kn_{y_i} = N(y_i)$ and $\kn_y = \km_y$ for all points $y$ of
codimension 1, distinct from all $y_i$.
\end{lemma}

{\bf Proof:}  One can easily construct $\kn$ with prescribed localizations as
in \cite{Bou} (VII 4.3).   Moreover, we may suppose that $\kn =
\cap_{codim\,y=1} \kn_y$.   But then $\kn$ is Cohen--Macaulay.

We can now prove the main result of this paper.   Recall that
\[
b(n,x) := b(n,
\kl(x)) = \max\{ \mbox{par}(n; \kl(x),\, \Gamma\}
\]
is the maximum number of independent parameters of isomorphism classes of
sandwiched $\kl(x)$--modules of rank $n$, which can be thought of as the
dimension of the ``moduli space'' of $\kl(x)$--CM--modules of rank $n$.

\begin{theorem}\label{4.9}
The function $b(n, x) = b(n, \kl(x))$ is upper semi--continuous.
\end{theorem}

{\bf Proof:}  Again we may suppose that $X$ is a smooth curve.   Let $g \in X$
be the generic point of $X$ and $\Lambda = \kl(g)$.   Find an overring $\Omega
\supset \Lambda$ such that $b(n, \Lambda) = {\rm par}(n, \Lambda, \Omega)$.
Using Lemma \ref{4.8}, we can construct a family of overrings $\kl' \supset
\kl$ with $\kl'(g) = \Omega$.   As $b(x) \ge {\rm par}(n, \kl(x), \kl'(x))$ for
every $x \in X$, it follows from Proposition \ref{4.7} that the set $\{x \in X
\mid b (x) \ge b(g)\}$ is closed.   This, of course, proves the theorem.

\begin{corollary}\label{4.10}
For any family of CM--algebras $(\kl, f : Y \to X),$ the set $W(\kl) = \{x
\in X \mid \kl(x)\mbox{ is wild}\}$ is a countable union of closed subsets of
$X$.
\end{corollary}

(For the definition of tame and wild CM--algebras cf.\ \cite{DG 1}).

The proof of this corollary follows from Theorem \ref{4.8} just in the same
way as it followed in the commutative case from Kn\"orrer's theorem (cf.\
\cite{DG 2}, Corollary 4.2).

Now we consider onother version of the semicontinuity theorem, where
algebras are given by parametrized families of generators. Namely,
let $X$ be an algebraic $k$--scheme, $\kl$ a family of CM--algebras with the
base $X$ and $\ki$ an ideal of $\kl$ such that $\kl/\ki$ is a locally free
$\ko_X$--module of finite rank, that is it corresponds to a vector bundle $\pi
: F \to X$.   The fibres $F(x)$ of $F$ are then finite--dimensional
$k(x)$--algebras.   Suppose given an algebraic $X$--scheme $f : Y \to X$ and
a set of $X$--morphisms $\{\gamma_i : Y \to F \mid i = 1, 2, \ldots, m\}$
(equivalently $Y$--sections of $f^\ast F$).

For each point $y \in Y$ denote $A(y)$ the subalgebra of $F(f(y))$ generated
by $\{\gamma_1(y),\; \gamma_2(y), \ldots, \gamma_m(y)\}$ and $\Lambda(y)$ the
preimage of $A(y)$ in $\kl(y) = \kl \otimes_{\ko_X} k(f(y))$.   Then
$\Lambda(y)$ is a CM--algebra, thus, given a family of overrings $\kl' \supset
\kl$, we may consider, as above, the functions on $Y$:
\vspace{-0.5cm}

\begin{eqnarray*}
p(n,d;\, y) &=& \mbox{par}(n,d;\, \Lambda(y),\; \kl'(y));\\
p(n;\, y) & = &\mbox{par}(n;\, \Lambda(y),\; \kl'(y));\\
b(n,y) & = &b(n, \Lambda(y)).
\end{eqnarray*}

\begin{theorem}\label{411}
In the above situation, the functions $p(n, d;\, y);\; p(n,y)$ and $b(n,y)$
are upper- semicontinuous on $Y$.
\end{theorem}

{\bf Proof}:  Replacing $\kl$ by $f^\ast(\kl)$, which is a family of
CM--algebras on $Y$, we may suppose that $X = Y$ and $f$ is the identity map.
Moreover, we may also suppose $X$ to be a smooth curve.   As the function
$\dim A(y)$ is obviously upper semi--continuous on $Y$, there is an open subset
$U \subset Y$ such that $\dim A(y)$ is constant and maximal possible on $U$.
Put $d = \dim F(y) - \dim A(y)$.   Then we obtain a section $\varphi : U \to
\mbox{Gr}(d, F)$ such that $A(y)$ is the subspace of $F(y)$ corresponding to
$\varphi(y)$ for each $y \in U$.   But as $X$ is a smooth curve and
$\mbox{Gr}(d,F)$ is projective over $X$, $\varphi$ can be prolonged to a
section $\bar{\varphi} : X \to \mbox{Gr}(d,F)$ (it follows, for example, from
\cite{Ha}, Proposition III.9.8).

Now $\bar{\varphi}$ gives rise to a subbundle $A' \subset F$ of constant
  codimension $d$.   Denote by $\Lambda'$ its preimage in $\kl$.   Note that
  both conditions
\[
\mbox{``}A'(x) \mbox{is a subalgebra of } F(x)\mbox{'' and ``}A'(x) \supset
A(x)\mbox{''}
\]
are evidently closed and hold on $U$.   Thus they hold on $X$, that is
$\Lambda'(x)$ is a subalgebra of $\kl(x)$ containing $\Lambda(x)$.   As
$\kl/\Lambda'$ is locally free of finite rank, $\Lambda'$ is really a family
of CM--algebras on $X$.   Hence, the functions $p'(n,d;\, x),\; p'(n;\, x)$
and $b'(n,x)$ defined just as $p(n,d;\, x),\; p(n;\, x)$ and $b(n,x)$ but
using $\Lambda'(x)$ instead of $\Lambda(x)$ are upper semicontinuous. On the
other hand we have inequalitites $p(n,d;\, x) \ge p'(n, d;\, x),\; p(n, x) \ge
p'(n, x),\; b(n,x) \ge b'(n,x)$
on $X$ and equality on $U$.   Therefore, $p(n,d;\, x),\; p(n;\, x)$
and $b(n;x)$ are also upper semicontinuous.

To show an application of Theorem~\ref{411}, we extend the criteria of
tameness, proved in \cite{DG 2} for the case char $k \not= 2$, to all
characteristics.
In order to do this, we must first define the singularities $T_{pq}$ in
positive characteristic, which are defined for char $k = 0$ as factorrings

$(\ast)\hspace{4cm}k[[X,Y]]/(X^p + Y^q + \lambda X^2 Y^2)$.

For our purpose  it is more convenient to define them using their
parametrization given by Schappert \cite{Sch}.   Namely, let $\Lambda$ be a
local commutative CM--algebra, $\Lambda_0$ its maximal overring.   Then
$\Lambda_0$ is a direct product of power series rings:
\[
\Lambda_0 \simeq k[[t_1]] \times k[[t_2]] \times \dots \times k[[t_s]]
\]
($s$ is ``the number of branches'' of $\Lambda$).   If $a \in \Lambda$, $a =
(a_1, a_2, \dots, a_s)$, with $a_i \in k[[t_i]]$, put {\boldmath$v$}$(a) =
(v(a_1), \dots, v(a_2))$, where $v(a_i)$ denotes the usual valuation on the
power series (in particular $v(0) = \infty$).

Call $\Lambda$ a plane curve singularity if its maximal ideal $\km$ is
generated by two elements:  $\km = (x,y)$.   Define the ({\bf valuation}) {\bf
  type} of $\Lambda$ as the pair ({\boldmath$v$}$(x)$, {\boldmath$v$}$(y)$).

\begin{definition}{\rm
Let $\Lambda$ be a plane curve singularity.   We say that $\Lambda$ is of {\bf
  type} {\boldmath$T_{pq}$}, where $p, q \in \N,\; \frac{1}{p} + \frac{1}{q}
  \le \frac{1}{2}$, if its valuation type is:
\[
\begin{array}{ll}
(2, p-2),\; (q-2, 2) & \mbox{for $p,q$ both odd},\\
(1, 1, p-2),\; (\infty, \frac{q}{2} - 1, 2) & \mbox{for $p$ odd, $q$ even},\\
(1, 1, \frac{p}{2} - 1, \infty), \; (\frac{q}{2} - 1, \infty, 1, 1) & \mbox{for
  $p,q$ both even}.
\end{array}
\]

By \cite{Sch}  this definition is equivalent to
the equation ($\ast$) if char $k = 0$.}
\end{definition}

The following theorem was proved in \cite{DG 2} for char $k \not= 2$.

Let $\Lambda$ be a local commutative CM--algebra, $\Lambda_0 = k[[t_1]] \times
\dots \times k[[t_s]]$ its maximal overring, {\boldmath$m$} = rad$\Lambda$.
Denote $t = (t_1, \dots, t_s) \in \Lambda_0$; $\Lambda^\prime = t\Lambda_0 +
\Lambda,\; \Lambda^{\prime\prime} = t${\boldmath$m$}$\Lambda_0 + \Lambda$
and $\Lambda^\prime_e = \Lambda^\prime + ke$, where $e \in \Lambda^\prime$ is
an idempotent.   For each overring $\Gamma \supset \Lambda^\prime$, let
$\Gamma/${\boldmath$m$}$\Gamma = L_1 \times \dots \times L_m$, where $L_i$ are
local algebras, $d_i = \dim\, L_i$, {\boldmath$d$}$(\Gamma) = (d_1, d_2, \dots,
d_m)$ and $d(\Gamma) = d_1 + \dots + d_m$ (the minimal number of
generators of $\Gamma$ as $\Lambda$--module).

\begin{theorem}
  If $\Lambda$ is of infinite Cohen--Macaulay type, the following conditions
  are equivalent:
\begin{enumerate}
\item[(1)] $\Lambda$ is tame.
\item[(2)] $\Lambda$ dominates a plane curve singularity of type $T_{pq}$ for
  some $p,q$ (that is, $\Lambda$ is isomorphic to an overring of $T_{pq}$).
\item[(3)] The following restrictions hold:
  \begin{enumerate}
  \item $d(\Lambda_0) \le 4$ and {\boldmath$d$}$(\Lambda_0) \not\in \{(4),\;
    (3,1),\; (3)\}$,
  \item $d(\Lambda^\prime) \le 3$ and {\boldmath$d$}$(\Lambda^\prime_e) \not=
    (3,1)$ for any idempotent $e$,
  \item if $d(\Lambda_0) = 3$, then $d(\Lambda^{\prime\prime}) \le 2$.
  \end{enumerate}
\end{enumerate}
\end{theorem}

{\bf Proof}: (1) $\Rightarrow$ (3) and (3) $\Rightarrow$ (2) were proved in
\cite{DG 2} and their proofs did not use the restriction char $k \not= 2$.  In
order to prove (2) $\Rightarrow$ (1), again following \cite{DG 2}, note that
the singularity $\Lambda$ of type $T_{pq}$ contains the $\Lambda_0$--ideal $I
= b\Lambda_0$, where
\[
\begin{array}{ll}
  b = (t_1^{p+1},\; t_2^{q+1}) & \mbox{for $p,q$ both odd};\\
  b = (t_1^{q/2+1},\; t_2^{q/2+1},\; t_3^{p+1}) & \mbox{for $p$ odd, $q$
  even};\\
  b = (t_1^{q/2+1},\; t_2^{p/2+1},\; t_3^{q/2+1},\; t_4^{p/2+1}) & \mbox{for
  $p,q$ both even}.
\end{array}
\]

Consider now a new CM--algebra $\Lambda(\lambda),\; \lambda \in k$,
containing $I$ and generated
modulo $I$ by the following 3 elements:
    \[
    \begin{array}{ll}
      (\lambda, 1)x,\; (1,\lambda)y,\; xy & \mbox{for $p,q$ both odd},\\
      (1,1,\lambda)x,\; (\lambda,\lambda, 1)y,\; xy & \mbox{for $p$ odd, $q$
        even},\\ (1, \lambda,1,\lambda)x,\; (\lambda,1, \lambda, 1)y,\, xy &
      \mbox{for $p,q$ both even}.
    \end{array}
    \]

If $(p,q) \not\in \{(4,4),\; (3,6)\}$, one can easily check that
$\Lambda(\lambda) \simeq \Lambda$ for $\lambda \not= 0$, while $\Lambda(0)$ is
a singularity of type $P_{pq}$ in the terminology of \cite{DG 2}, that
is generated modulo $\,I \,$ by the elements $\,x_0,y_o \,$ such that
{\boldmath$v$}$(x_0)$,{\boldmath$v$}$(y_0)$ are of the form:
\[
\begin{array}{ll}
(2, \infty),\; (\infty, 2) & \mbox{for $p,q$ both odd},\\
(1, 1, \infty),\; (\infty, \infty, 2) & \mbox{for $p$ odd, $q$ even},\\
(1, 1, \infty, \infty), \; (\infty, \infty, 1, 1) & \mbox{for $p,q$ both even}.
\end{array}
\]

   Again, the
calculations for $P_{pq}$ in \cite{DG 2} did not use the condition char $k\not=
2$.   Hence, they are tame and Theorem 4.11 implies that $\Lambda$ is also
tame.   The calculation of Dieterich for the remaining case $(p,q) = (3,6)$ or
$(p,q) = (4,4)$ (cf.\ \cite{Di 1}, \cite{Di 2}) also did not use any
conditions on characteristics.  Thus, implication (2) $\Rightarrow$ (1) is
completely proved.

\newpage

\section{Some analogues}

Here we give some  examples of  ``almost versal families'' and
 semicontinuity theorems for
other situations in representation theory. As all the proofs are quite
similar (and easier) to those of the preceding sections,
we omit them and give only the final formulations of the results
analogous to Theorems~\ref{3.4a}, \ref{4.7} and
\ref{411}. Although some of the corresponding semicontinuity theorems
are known, we hope that the ``unification'' will be useful for these
cases too. At least we give  new proofs for them.

\subsection*{Finite-dimensional algebras}

Here, let  $A$ be a finite--dimensional algebra over an algebraically
closed field $k$. A {\bf family of {\boldmath$A$}-modules\/} parametrized
by a $ k $--scheme $ X  $ is a sheaf $ \km  $ of
$ A\otimes_k\ko_X$--modules, which is coherent and flat over $ \ko_X
 $.   To be in the frame of projective varieties, we can consider first
the subvariety $B(P, I, d) \subset Gr(P, d)$, where $P$ is a
projective module over a finite--dimensional algebra $A$, $I$ an ideal
of $A$ contained in the radical and $B(P, I, d)$ consists of the
$A$--submodules $L \subset IP$. $Gr(P,d)$ denotes the Grassmannian of
$d$--codimensional subspaces of $P$.   Then the canonical sheaf
$ \kf=\kf(P,I,d)  $ on $ B(P,I,d)  $ is a family of $ A $--modules
and the following result holds.

\begin{theorem}
\label{5.1}
Let $ A  $ be a finite-dimensional $ k $-algebra, $ J=rad\, A  $ and
 $ \km  $  a family of $ A $--modules parametrized by a reduced
$ k $--scheme $ X  $. Then there exists a descending chain of closed
subschemes $ X=X_0\supset X_1\supset X_2\supset \dots\supset
X_m=\emptyset  $ and a set of morphisms $ \{\varphi_i:Y_i\rightarrow
B(P_i,J,d_i)\mid  i=1,\dots,m\}  $ for some projective $ A $--modules $ P_i  $
 such that $ \km_{Y_i}\simeq\varphi_i^\ast \kf(P_i,J,d_i)  $, where
$ Y_i=X_{i-1}\setminus X_i  $. Moreover, if $ r= \mbox{ rank } \km  $ (as a
locally free sheaf over $ X  $), then $ \dim P_i\le rp  $, where
$ p  $ is the maximal dimension of indecomposable projective
$ A $--modules.
\end{theorem}

The group $ G = \mbox{ Aut}_AP  $ acts on $ B=B(P,I,d)  $ and, as $ I\subset
\mbox{ rad } A  $, we conclude that $ \kf(x)\simeq\kf(y)  $ if and only
if $ x  $ and $ y  $ belong to the same $ G $--orbit. The subsets
$ B_i=\{ x\in B  |  \dim(Gx)\le i \}  $ are obviously closed in $ B
 $. Hence, we can define the {\bf number of parameters\/}:
\[
{\rm par}(P,I,d; A) = \max_i(\dim B_i-i)
\]
and
\[
{\rm par}(P,I; A) = \max_d {\rm par}(P,I,d;,A).
\]
In particular, put
\[
{\rm par}(n,d; A)={\rm par}(nA,rad A,d; A)\quad {\rm and}\quad
{\rm par}(n; A)={\rm par}(nA,rad A; A).
\]

Just as for Cohen-Macaulay algebras, these numbers  give upper bounds
for the number of (independent) parameters of isomorphism classes of
$A$--modules of rank $n$ in
{\em any\/} family of $ A $--modules.

Consider  now a {\bf family of algebras} parametrized by a $ k
$--scheme $ X  $, that is a flat coherent sheaf of $ \ko_X$-algebras $ \ka
$. Then we are able to define the following functions on
$ X  $:
\vspace{-0.5cm}

\begin{eqnarray*}
        {\rm par}(x,n,d) &=& {\rm par} (n,d; \ka(x))  ,\\
        {\rm par}(x,n)   &=& {\rm par} (n; \ka(x)) .
\end{eqnarray*}

\begin{theorem}
\label{5.2}
For each family of finite--dimensional $ k $--algebras the functions
$ {\rm par}(x,n,d)  $ and $ {\rm par}(x,n)  $ are upper--semicontinuous.
\end{theorem}

A version of this theorem was proved by Gei{\ss} \cite{Gei}.   Note also that
Theorem 5.2 provides a new proof of Gabriel's theorem \cite{Gab} that finite
representation type is an open condition.   This follows since the
Brauer--Thrall conjectures are known to be true for finite dimensional
algebras.

It is easy  to generalize the last theorem to the situation
where the algebras are given by ``generators and relations''.  Namely,
suppose we are given:
\begin{itemize}
\item
a family $ \ka  $ of finite--dimensional $ k $-algebras over $X$;
\item
an algebraic $ X $-scheme $ f:Y\rightarrow X   $;
\item
two sets of $ X $-morphisms $ \{\gamma_i:Y\rightarrow
F \mid  i=1,\dots,m\} $ and $ \{\rho_j:Y\rightarrow
F \mid  j=1,\dots,r\} $, where $ F  $ is the  vector bundle on $ X
 $ corresponding to the locally free sheaf $ \ka  $.
\end{itemize}

For any point $ y\in Y  $, denote by $ I(y)  $ the ideal in $ F(f(y))
 $ generated by the set $ \{\rho_j(y) | j=1,\dots,r \} $ and by $ A(y)
 $ the subalgebra of $ F(f(y))/I(y)  $ generated by the classes
$ \{\gamma_i(y)+I(y) | i=1,\dots,m \} $. Then we can define the
functions on $ Y  $:
\vspace{-0.5cm}

\begin{eqnarray*}
        p(y,n,d) &=& {\rm par} (n,d; A(y))  ,\\
        p(y,n)   &=& {\rm par} (n; A(y)) .
\end{eqnarray*}

\begin{theorem}
\label{5.3}
In the above situation the functions $ p(y,n,d)  $ and $ p(y,n)  $
are upper-semicontinuous on $ Y  $.
\end{theorem}

\subsection*{Bimodules}

Consider now the categories of elements of
  finite--dimensional bimodules (in the sense of \cite{D2}, although we
give here a somewhat different definition).  Let $ A
   $  be a finite--dimensional $ k $--algebra, where
  $ k  $ is again an algebraically closed field, and let $ V  $ be a
  finite--dimensional $ A $--bimodule. The {\bf
  elements} of $ V  $ are, by definition, those of the set
  $ El(V)=\bigsqcup_{P} V(P)  $, where $ P  $ runs
 through all (finitely generated) projective $ A $--modules
  and $ V(P)=\mbox{ Hom}_A(P,V\otimes_AP)  $. Two
  elements $ u\in V(P)  $ and $ u'\in V(P')  $ are said to be
  {\bf isomorphic} if there exists an isomorphism $ p:P\rightarrow P' $
  such that $ u'=(1\otimes p)up^{-1} $. Indeed, in \cite{D2} only so--called
  disjoint bimodules were considered. The bimodule $ V $ is
said to be {\bf disjoint} if $ A=A_1\times A_2  $ and $ VA_1=A_2V=0
 $. In most applications in  representations theory one needs only
disjoint bimodules, but non--disjoint ones appear in various
``reduction processes''.

 To remain in the category of
projective varieties, it is convenient to change  the
problem slightly. Namely, call elements $ u  $ and $ u'  $ {\bf
(projectively)  equivalent}, if $ u'  $ is isomorphic to $ \lambda u
 $ for some non--zero $ \lambda\in k  $. Obviously, if the bimodule
is disjoint, then equivalent elements are isomorphic, but in the
non--disjoint case it is not always so.

Let $ X  $ be a $ k $--scheme and $ \kp  $ a flat coherent sheaf
of $ A_X $--modules, where $ A_X=A\otimes_k\ko_X  $. Put
$ V(\kp)=\kh{\it om}_{A_X}(\kp,V\otimes_A\kp)  $. It is a locally free
coherent sheaf of $ \ko_X$--modules. Hence, the corresponding
projective bundle $ \P_\kp=\P_X(V(\kp))  $ over $ X  $ is defined (cf.\
\cite{Ha}).
A {\bf projective family} (or simply {\bf family}, as we do not
consider other families here) of elements of the bimodule $ V  $
with  base $ X  $ is, by definition, a section
$ \Phi:X\rightarrow \P_\kp $ for some $ \kp  $.  Note that
using  projective families we need to consider projective
equivalence instead of isomorphism and to exclude zero elements of the
bimodule. But this does not essentially differ from the
classification problem for the elements of bimodules up to isomorphism.

The ``almost universal'' families in this case are more or less
evident. Indeed, put, for any projective $ A $--module $ P  $,
$ B(P)=\P_k(V(P))  $ and $ \tilde P=P\otimes_k\ko_{B(P)}  $. Then
$ \P_{\tilde P}\simeq B\times B  $, where $ B=B(P)  $ and the
diagonal map $ \Delta_P:B\rightarrow B\times B $ defines a family of
elements of $ V  $ with the base $ B  $. The following result is
almost obvious.

\begin{theorem}\label{5.4}
Let $ A  $ be a finite--dimensional $ k $--algebra, $ V  $ a
finite-dimensional $ A$--bimodule and and $ \Phi:X\rightarrow
\P_\kp $ a (projective) family of elements of $ V  $.
 Then there exists a descending chain of closed
subschemes $ X=X_0\supset X_1\supset X_2\supset \dots\supset
X_m=\emptyset  $ and a set of morphisms $ \{\varphi_i:Y_i\rightarrow
B(P_i) | i=1,\dots,m\}  $ for some projective $ A $--modules $ P_i  $ such
that $ \kp_{Y_i}\simeq P_i\otimes_k\ko_{Y_i}$.   Hence, the
restriction of $ \P_\kp  $ on $ Y_i= X_{i-1}\backslash X_i $ can be identified
with $ Y_i\times B(P_i)  $, and, under this identification,
$ \Phi_{Y_i}=1\times\varphi_i  $. Moreover, $ \dim P_i=rank \kp  $
for all $ i  $.
\end{theorem}

The group $ G=\mbox{ Aut}_AP  $ acts on $ B=B(P)  $ and its orbits are the
classes of projective equivalence. Hence, we are again able to define the
closed subsets
$ B_i=\{ x\in B  |  \dim(Gx)\le i  \} $ and the {\bf number of parameters\/}:
\[
{\rm par}(P; A,V) = \max_i(\dim B_i-i) ,
\]
in particular
\[
{\rm par}(n; A,V)=par(nA; A,V) .
\]

Now, given a family of algebras $ \ka  $ with base  $ X  $ and
a family of bimodules, that is a coherent sheaf $ \kv  $ of $ \ka
$--bimodules, flat over $ \ko_X  $,
we can define the function on
$ X  $:
\[
{\rm  par}(x,n) =   {\rm par} (n; \ka(x),\kv(x)) .
\]

\begin{theorem}
\label{5.5}
For each family of finite--dimensional $ k $--algebras and bimodules
the function
$ {\rm par}(x,n)  $ is upper--semicontinuous.
\end{theorem}

Of course, one could easily give a version of the last theorem,
where the algebras and bimodules are defined by generators and
relations, but we leave this obvious generalization to the reader.

\subsection*{Remark}

 In particular, in both cases we can see that the set of wild algebras
  (or bimodules) in some family is again a countable union of closed
  subsets.  It looks very likely that this set is even closed and,
  hence, that the set of tame algebras (or bimodules) is open.  In
  order to prove it, one only needs to show that the set of tame
  algebras (bimodules) is really a countable union of constructible
  sets (cf.\  \cite{Gab}).

 If we consider families of commutative CM--algebras, then the set of tame
  algebras is indeed open.  This can be derived from \cite{DG 2} in two ways.
  The first is to apply the classification of \cite{DG 2} and deformation
  theory of singularities:   the set of singularities which are of finite
  CM--representation type or which are tame is open in any flat family of
  singularities.   The second is to note that the strict
  respresentations over the free algebra $k\langle x,y\rangle$ constructed in
  \cite{DG 2} are of
  bounded rank.  Hence, we can find a common constant $n$ such that a
  commutative CM--algebra $\Lambda$ is wild if and only if $p(n,\Lambda) >
  rn$, where $r$ is the rational length of $\Lambda$, which coincides in this
  case with the number of branches.  As $r$ is obviously bounded in any
  family, we have now only to apply Theorem 4.9.

\end{document}